\documentclass[11 pt]{amsart}
\usepackage{amsmath,amsthm,amsfonts,amssymb,amscd,amsrefs}
\usepackage[top=3cm, bottom=2.9cm,left=2cm, right=2cm]{geometry}
\allowdisplaybreaks[1]%this allows to break lines in \align .. 
\usepackage{centernot}
\usepackage{xcolor} %not originally present added to higlight the text
\usepackage{soul} %not originally present added to higlight the text
\usepackage[linkcolor=red,citecolor=red,colorlinks=true]{hyperref}
\usepackage{mathrsfs}
\usepackage{mathtools}
\usepackage{etoolbox}% for the patchcmd command i.e. to get rid of d small caps style in title
\usepackage{scalerel,stackengine}
\stackMath
\newcommand\widecheck[1]{%
	\savestack{\tmpbox}{\stretchto{%
			\scaleto{%
				\scalerel*[\widthof{\ensuremath{#1}}]{\kern-.6pt\bigwedge\kern-.6pt}%
				{\rule[-\textheight/2]{1ex}{\textheight}}%WIDTH-LIMITED BIG WEDGE
			}{\textheight}% 
			
		}{0.5ex}}%
	\stackon[1pt]{#1}{\scalebox{-1}{\tmpbox}}%
}
\theoremstyle{definition}
\newtheorem{defn}{Definition}[section]

\newtheorem{eg}[defn]{Example}

\newtheorem{case}{Case}%will come as case1 case2,...; put * for no number; use [defn] to give no, within section.
\theoremstyle{remark}

\newtheorem{rmk}[defn]{Remark}
\theoremstyle{plain}
\newtheorem{thm}[defn]{Theorem}

\newtheorem{cor}[defn]{Corollary}
\newtheorem{prop}[defn]{Proposition}
\DeclareMathOperator*{\spn}{span}

\numberwithin{equation}{section}
\makeatletter %for roman letter

\newcommand{\Rmnum}[1]{\expandafter\@slowromancap\romannumeral #1@}

\makeatother
\makeatletter
\def\imod#1{\allowbreak\mkern10mu({\operator@font mod}\,\,#1)}
\makeatother
\makeatletter % for making title to appear as we hav written(otherwise it is small caps)
\patchcmd{\@settitle}{\uppercasenonmath\@title}{}{}{}
\patchcmd{\@setauthors}{\MakeUppercase}{}{}{}
\patchcmd{\section}{\scshape}{}{}{}
\makeatother
\begin{document}
	% \title[Bhatia-\v Semrl Property, Strong Subdifferentiability and Essential Norm]{On Bhatia-\v Semrl Property, Strong Subdifferentiability and Essential Norm of Operators on Banach Spaces}
    % \title[]{On Essential Norm, Birkhoff-James Orthogonality and Strong Subdifferentiability   in Space of  Operators on Banach Spaces}
     \title[Essential Norm, Bhatia-\v{S}emrl Property and Strong Subdifferentiability]{Connections between Essential Norm, Bhatia-\v{S}emrl Property and Strong Subdifferentiability of Operators on Banach Spaces}
	\author{C. R. Jayanarayanan}
    \author{Rishit R Rajpopat}
	\address{Department of Mathematics, Indian Institute of Technology Palakkad, 678557, India}
	%\curraddr{Department of Mathematics, Indian Institute of Technology Palakkad, 678557, India}
	\email{crjayan@iitpkd.ac.in}
	% \address{Department of Mathematics, Indian Institute of Technology Palakkad, 678557, India}
	\email{rishitraj3@gmail.com, 212214004@smail.iitpkd.ac.in}
	\begin{abstract}
    % \fontsize{9.5}{11}\selectfont
    { We study  Birkhoff-James orthogonality in  spaces of bounded linear operators between Banach spaces through  the Bhatia-\v Semrl property.  We investigate the interplay among three key properties of bounded linear operators: strong subdifferentiability, the Bhatia-\v Semrl property,  and  the condition that the essential norm is strictly less than the operator norm. 
    % For a Hilbert space $H$ and for $1<p,q<\infty$, we show that for any nonzero operator in $B(H)$ and $B(\ell^p, \ell^q)$,  the essential norm is strictly less than the operator norm if and only if it is a point of strong subdifferentiability of the norm and its norm-attainment set is compact. 
     If  $X$ is a Banach space with a $1$-unconditional basis and  $Y$ is any Banach space, we show that every nonzero operator  $T \in B(X, Y)$ satisfying the Bhatia-\v{S}emrl property satisfies $\|T\|_e<\|T\|$. We identify a class of pairs of Banach spaces  $X$ and $Y$ for which the inequality $\|T\|_e<\|T\|$ can be characterized in terms of strong subdifferentiability of the operator norm at $T$ and compactness of the norm attainment set of $T$.
    % Moreover, for nonzero operators in these spaces that satisfy the Bhatia-\v Semrl property, we show that their essential norm must  be strictly less than their operator norm. 
    We also study norm one projections satisfying the Bhatia-\v Semrl property and provide examples of operators that possess this property.}
	\end{abstract}
	\maketitle
\section{Introduction}
The notion of orthogonality plays an important role in understanding the geometric properties of Hilbert spaces. Since a Banach space is not necessarily an inner product space, there is no canonical way to define orthogonality in such spaces. Several notions of orthogonality in Banach spaces have been introduced in the literature  and among them, Birkhoff-James orthogonality, introduced by Birkhoff in \cite{B}, has emerged as an important tool in understanding the geometry of Banach spaces due to its relation with several geometric properties such as  strict convexity, smoothness, differentiability etc. In this article, we investigate the geometric properties of elements of the space of all bounded linear operators between Banach spaces through the framework of Birkhoff-James orthogonality.

In this article, we consider only Banach spaces over the real field $\mathbb{R}$. We consider every Banach space $X$, under the canonical embedding, as a subspace of $X^{**}$. For a Banach space $X$; let $B_X$ and $S_X$ denote the closed unit ball and  the unit sphere of $X$ respectively. For Banach spaces $X$ and $Y$, let ${B}(X, Y)$ denote the space of all bounded linear operators from $X$ to $Y$, and ${K}(X, Y)$ denote the space of all compact operators from $X$ to $Y$. When $X=Y$, we denote $B(X, Y)$ by $B(X)$ and $K(X, Y)$ by $K(X)$. For $T\in B(X, Y)$, let $M_T$ denote the set of all unit vectors in $X$ at which $T$ attains its norm, that is, $M_T=\{x\in S_X: \|Tx\|=\|T\|\}$. An operator $T \in B(X, Y)$ is said to be norm-attaining if $M_T \neq \emptyset$.   The essential norm $\|T\|_e$ of a bounded operator $T\in {B}(X, Y)$ is the distance from $T$ to $K(X, Y)$, that is, $\|T\|_e = d(T, {K}(X, Y))$.

Let $X$ be a Banach space and $x,y \in X$. We say that $x$ is  Birkhoff-James orthogonal to $y$, denoted by $x \perp_B y$, if  $\|x+\lambda y\| \geq \|x\|$  for every scalar $\lambda$ (see \cite{B}). It is easy to observe that in an inner product space, the Birkhoff-James orthogonality coincides with the usual inner product orthogonality. Birkhoff-James orthogonality in the space of bounded linear operators between Banach spaces has been extensively studied in the literature (see \cites{MPSbook, K, KL, SPH}). In \cite{BS}, Bhatia and \v Semrl proved that for a Hilbert space $H$  and for operators $T, A\in B(H)$, $T \perp_B A$ if and only if there exists a sequence $(x_n)$ of unit vectors in $H$ such that $\|Tx_n\|\to \|T\|$ and $\langle Tx_n, Ax_n \rangle \to 0$. Furthermore, they have also proved in \cite{BS} that if $H$ is finite dimensional, then $T \perp_B A$ if and only if there exists a vector $x\in M_T$ such that $\langle Tx, Ax \rangle=0$. This result of Bhatia and \v{S}emrl motivated the introduction of the following notion in \cite{SPH}.
\begin{defn}
Let $X$ and $Y$ be Banach spaces. An operator $T \in B(X, Y)$ is said to satisfy the Bhatia-\v Semrl property if for every $A \in B(X, Y)$ with $T \perp_B A$, there exists $x \in M_T$ such that $Tx \perp_B Ax$.
\end{defn}

 In \cite{BS}, Bhatia and \v Semrl conjectured that if $X$ is  the Banach space $\mathbb{C}^n$ equipped with an arbitrary  norm, then every operator on $X$ satisfies the Bhatia-\v Semrl property. However, in \cite{LS}, Li and Schneider  provided counterexamples to show that this conjecture is not valid in general.  In fact, Benítez et al. proved in \cite{BFS} that a finite dimensional real normed linear space $X$ is an inner product space if and only if every $T \in B(X)$ satisfies the Bhatia-\v Semrl property. Moreover, it was observed in \cite{KL} that a Banach space $X$ is reflexive if and only if every $x^* \in X^*$ satisfies the Bhatia-\v{S}emrl property. 

The following result characterizes operators on infinite dimensional Hilbert spaces which satisfy the Bhatia-\v Semrl property.

\begin{thm}[\cite{PSG}*{Theorem 3.1}]
\label{bshs}
    Let $H$ be an infinite dimensional Hilbert space and $T \in B(H)$ be such that $T \neq 0$. Then $T$ satisfies the Bhatia-\v Semrl property if and only if there exists a finite dimensional subspace $H_0$ of $H$ such that $M_T=S_{H_{0}}$ and $\|T | _{{H_0}^{\perp}}\| < \|T\|$.
\end{thm}
In \cite{CK}*{Proposition 3.7}, it was shown that there does not exist any nonzero $T \in B(c_0)$ that satisfies the Bhatia-\v Semrl property. Similarly, it was shown in \cite{KL}*{Proposition 3.6}  that for any Banach space $Y$,  there does not exist any nonzero $T\in B(L^1[0, 1], Y)$ that satisfies the Bhatia-\v Semrl property. However,  if $X$ has the Radon-Nikod\' ym property, then it was shown in \cite{K}*{Theorem 2.4} that the set of all norm attaining operators in $B(X, Y)$ that satisfies the  Bhatia-\v Semrl property forms a dense set. 

For a closed subspace $J$ of a Banach space $X$ and $x\in X$, we say that $x$ is Birkhoff-James orthogonal to $J$, denoted by $x \perp_B J$, if $x \perp_B y $ for every $y\in J$. One can easily observe that  $x \perp_B J$ if and only $d(x, J)=\|x\|$. In particular, for Banach spaces $X$ and $Y$, and for any $T\in B(X,Y)$, we have $T\perp_B K(X, Y)$ if and only if $\|T\|_e=\|T\|$. 
% Motivated by these observations, we investigate the relationship between the condition $\|T\|_e=\|T\|$ and the Bhatia-\v Semrl property of $T$. 

% Another pair of Banach spaces $(X, Y)$ for which this property is partially understood  is when $K(X, Y)$ is an $M$-ideal in $B(X, Y)$.
The next result from \cite{MPRS} provides a connection between the condition $\|T\|_e<\|T\|$ and the Bhatia-\v Semrl property of $T$ when $K(X, Y)$ is an $M$-ideal in $B(X, Y)$. We first recall the notion of an $M$-ideal in a Banach space. 
% \begin{defn}[\cite{hww}]
% Let $X$ be a Banach space and $Y$ be a closed subspace of $X$. A linear projection P on $X$ is said to be an \emph{$L$-projection} if  $\|x\|=\|Px\|+\|x-Px\|$ for all $x \in X$. A subspace $Y$ of $X$ is said to be an \emph{$M$-ideal}  in $X$ if $Y^\bot$ is the range of an $L$-projection in $X^*$. 
% \end{defn}
% \begin{defn}[\cite{hww}]
% Let $X$ be a Banach space and let  $P: X\to X$ be a linear projection on $X$. Then $P$ is said to be an $L$-projection ($M$-projection) if $\|x\|=\|Px\|+\|x-Px\|$
%   ($\|x\|=\max\{\|Px\|,\|x-Px\|\}$) for all $x \in X$. For $1\leq p <\infty$, we say $P$ is an $L^p$-projection if $\|x\|^p = \|Px\|^p+\|x-Px\|^p$ for all $x \in X$. A closed subspace $Y$ of $X$ is said to be an $L$-summand in $X$  if it is the range of an $L$-projection. A closed subspace $Y$ of $X$ is said to be an $L^p$-summand in $X$ if it is the range of an $L^p$-projection. A closed subspace $Y$ of $X$ is said to be an {$M$-ideal} in $X$ if $Y^\bot$ is an $L$-summand  in $X^*$.
%   %A Banach space $X$ is said to be $M$-embedded space if $X$ is an $M$-ideal in $X^{**}$.
% \end{defn}

\begin{defn}[\cites{hww, BDE}]
Let $X$ be a Banach space and let  $P\in B(X)$ be a projection.  For $1\leq p <\infty$, we say $P$ is an $L^p$-projection ($M$-projection) if $\|x\|^p = \|Px\|^p+\|x-Px\|^p$ ($\|x\|=\max\{\|Px\|,\|x-Px\|\}$) for all $x \in X$. 
% A closed subspace $Y$ of $X$ is said to be an $L^p$-summand ($M$-summand) in $X$ if it is the range of an $L^p$-projection ($M$-projection).
We call an $L^1$-projection as an $L$-projection.
% and an $L^1$-summand as an $L$-summand. 
A closed subspace $Y$ of a Banach space $X$ is said to be an {$M$-ideal} in $X$ if there exists an  $L$-projection $P : X^* \to X^*$ such that the range of $P$ is $Y^{\perp}$ (the annihilator of $Y$ in $X^*$).
  %A Banach space $X$ is said to be $M$-embedded space if $X$ is an $M$-ideal in $X^{**}$.
\end{defn}

% The following theorem from \cite{MPRS} provides a connection between the condition $\|T\|_e<\|T\|$ and the Bhatia-\v Semrl property of $T$ when $K(X, Y)$ is an $M$-ideal in $B(X, Y)$.
\begin{thm}[\cite{MPRS}*{Theorem 3.2}]
    Let $X$ be a reflexive Banach space and $Y$ be any Banach space, such that $K(X, Y)$ is an $M$-ideal in $B(X, Y)$. Let $T \in S_{B(X, Y)}$ be such that $\|T\|_e <1$ and $M_T= D \cup (-D)$, where $D$ is a closed connected subset of $S_X$. Then $T$ satisfies the Bhatia-\v Semrl property.
\end{thm}
% Naturally one can ask, given any $X, Y$ whether there exists $T \neq 0 \in B(X, Y)$  that satisfies Bhatia-\v Semrl property $?$  Unfortunately the answer to this negative.

One of the motivations for studying Birkhoff-James  orthogonality is its deep connection with the differentiability properties of the norm (see \cite{J}*{Theorem 4.2}). 

We say  that the norm of a Banach space $X$ is G\^ateaux differentiable at a vector $x \in X$ if  
$
\lim _{t \rightarrow 0} \frac{\|x+t y\|-\|x\|}{t} 
$
exists for every $y \in X$. The norm of $X$ is said to be  Fr\'echet differentiable at   $x\in X$ if this limit is uniform for every $y \in B_X$. The norm of $X$ is called Fr\'echet (resp. G\^ateaux) differentiable if it is Fr\'echet (resp. G\^ateaux) differentiable at every point of $X \backslash\{0\}$. 
A vector $x \in X$ is said to be a smooth point of $X$ if there exists a unique $x^* \in S_{X^*}$ such that $x^*(x)=\|x\|$. It is well known that the norm of $X$ is G\^ateaux differentiable at $x\in X$ if and only if $x$ is a smooth point of $X$. The Banach space $X$ is said to be a smooth space if every element of $S_{X}$ is smooth. Similarly, a Banach space $X$ is said to be uniformly smooth if $
 \lim _{t \rightarrow 0} \frac{\|x+t y\|-\|x\|}{t} 
 $ exists for every $x \in S_X$ and $y \in X$, and the convergence is uniform for $(x, y) \in S_X \times S_X$. We say that a Banach space $X$ is strictly convex if $\|\frac{x+y}{2}\| <1$ whenever $x$ and $y$ are distinct elements of $S_X$. Moreover, a Banach space $X$ is said to be uniformly convex if for every $0<\epsilon \leq 2$, there exists a $\delta >0$ such that $\left\|\frac{x+y}{2} \right\|\leq 1-\delta$ whenever $x, y \in S_X$ and $\|x-y\| \geq \epsilon$.

The following result characterizes the smooth points in $B(X, Y)$ using the Bhatia-\v Semrl property.
\begin{thm}[\cite{SPMR}*{Theorem 3.3}]
    Let $X$ and $Y$ be Banach spaces and let $T \in B(X, Y)$ be such that $T \neq 0$ and $M_T \neq \emptyset$. Then $T$ is a smooth point if and only if the following conditions hold:
    \begin{enumerate}
        \item $M_T= \{\pm x_0 \}$, for some $x_0 \in S_X$.
        \item $Tx_0$ is a smooth point in $Y$.
        \item $T$ satisfies the Bhatia-\v Semrl property. 
    \end{enumerate}
\end{thm}
The following non-smooth extension of Fr\'echet differentiability was introduced in \cite{FP}.
\begin{defn}[\cite{FP}]
 The norm of a Banach space $X$ is said to be  strongly subdifferentiable  at $x\in X$ 
if the one sided limit,
$
\lim_{t\to 0^+}\frac{\|x+ty\|-\|x\|}{t}
$
exists uniformly for $y\in B_X$. 
\end{defn}

In \cite{FP}, it is observed that the norm  of a Banach space $X$ is Fr\'echet differentiable at $x \in X$ if and only if it is both G\^ateaux differentiable and strongly subdifferentiable at $x$. In this article, we investigate the relationship among the condition $\|T\|_e<\|T\|$, the strong subdifferentiability, and the Bhatia-\v Semrl property of $T$.

It implicitly follows from Theorem~\ref{bshs} and  \cite{S} that, for a Hilbert space $H$, an operator $T \in B(H)$ satisfies the Bhatia-\v{S}emrl property if and only if $\|T\|_e<\|T\|$. However, for general Banach spaces, $\|T\|_e<\|T\|$ need not imply the Bhatia-\v{S}emrl property, as shown by Example~4.3 in \cite{LS}. In Section \ref{bsssd}, we prove that if  $X$ is a Banach space with a $1$-unconditional basis and  $Y$ is any Banach space, then every operator  $T \in B(X, Y)$ satisfying the Bhatia-\v{S}emrl property satisfies $\|T\|_e<\|T\|$. Example \ref{ceg} in Section \ref{bsssd} shows that the Bhatia-\v{S}emrl property alone does not guarantee the compactness of the norm attainment set. In Section \ref{bsssd}, we identify conditions on Banach spaces $X$ and $Y$ under which every operator with the  Bhatia-\v{S}emrl property has a compact norm attainment set. We show that if $X$ is a reflexive Banach space with the Kadets-Klee property and $Y$ is a Banach space such that $K(X, Y)$ is an $M$-ideal in $B(X, Y)$,  then $M_T$ is  non-empty and compact for every operator $T \in B(X, Y)$  with $\|T\|_e<\|T\|$. In fact, under some additional assumptions on  $X$ and $Y$, we characterize operators $T \in B(X, Y)$ satisfying the condition $\|T\|_e<\|T\|$ in terms of strongly subdifferentiability of the operator norm at $T$ and compactness of $M_T$.

% This article has three main sections. Before commenting on the contents of each section let us make some important observations. 
In Section \ref{bjp}, we study the Bhatia-\v{S}emrl property of norm one projections and isometries between Banach spaces. For an orthogonal projection $P$  on an infinite dimensional Hilbert space $H$, we have that $M_P=S_{R(P)}$, where $R(P)$ denotes the range of $P$. Then it follows from Theorem~\ref{bshs} that an orthogonal projection $P$ on an infinite dimensional Hilbert space $H$ satisfies the Bhatia-\v Semrl property if and only if $P$ is a finite rank projection. Similarly, we can show that  no isometry $T$ between an infinite dimensional Hilbert space $H$ satisfies the Bhatia-\v Semrl property as $M_T=S_H$. However, it is clear from \cite{CK}*{Proposition 3.7} and \cite{KL}*{Proposition 3.6} that finite rank norm one projections on $c_0$ and $L^1[0, 1]$ do not satisfy the Bhatia-\v Semrl property, respectively.  We prove that if $X$ is a reflexive, strictly convex Banach space with the Kadets-Klee property, then any finite rank norm one projection $P$ on $X$ satisfies the Bhatia-\v{S}emrl property. We also obtain a similar result for finite rank $L^p$ projections $(1\leq p<\infty)$ without any assumptions on the Banach space $X$. Finally, we establish that if $\mu$ is a $\sigma$-finite positive measure and $L^p(\mu)$ is infinite dimensional, then no isometry from $L^p(\mu)$ to $L^p(\mu)$ satisfies the Bhatia-\v{S}emrl property.

\section{Bhatia-\v Semrl property, essential norm and strong subdifferentiability in space of bounded linear operators on Banach spaces}
\label{bsssd}
In this section, we investigate the relationship among the properties: the  Bhatia-\v Semrl property, the condition that the essential norm is less than the operator norm, and the strong subdifferentiability of operators between Banach spaces. It implicitly follows from Theorem~\ref{bshs} and  \cite{S} that for a Hilbert space $H$, an operator $T \in B(H)$ satisfies the Bhatia-\v{S}emrl property if and only if $\|T\|_e<\|T\|$.
The following result provides an explicit connection among these three properties in the context of Hilbert spaces. 
% The following result answer the Question \ref{q3} positively for $B(H)$ 

\begin{prop}\label{BSHESS}
    Let $H$ be a Hilbert space and $T \in S_{B(H)}$. Then the following are equivalent.
    \begin{enumerate}
        \item $T$ satisfies the Bhatia-\v Semrl property.
        \item  $\|T\|_e<1$.
        \item $M_T$ is compact and the norm of $B(H)$ is strongly subdifferentiable at $T$.
    \end{enumerate}
\end{prop}
\begin{proof} (1)$\implies$(2): Suppose $T$ satisfies the Bhatia-\v Semrl property. Then, by \cite{PSG}*{Theorem 3.1}, there exists a finite dimensional subspace $H_0$ such that $M_T=S_{H{_0}}$ and $\|T|_{{H_0}^{\perp}}\|<1$. Let $P$ be the orthogonal projection onto $H_0$. Then $\|T(I-P)\|=\|T |_{{H_0}^\perp}\|<1$. Since $P$ is of finite rank, we have that $TP$ is a compact operator. Therefore, $\|T\|_e \leq \|T-TP\|<1$.
   
(2)$\implies$(3): Suppose $\|T\|_e<1$. Then, by \cite{S}*{Corollary 2.6}, it follows that $T$ satisfies the Bhatia-\v Semrl property. Hence, by \cite{PSG}*{Theorem 3.1}, there exists a finite dimensional subspace $H_0$ such that $M_T=S_{H{_0}}$ and $\|T|_{{H_0}^{\perp}}\|<1$. Since $H_0$ is finite dimensional, we have that $M_T$ is compact.  Now let $|T|=(T^* T)^\frac{1}{2}$ and let $x\in S_{H_0^\perp}$. Then $ \||T|(x)\|=\|Tx\|$ and hence $\|\,|T|\big | _{{H_{0}}^\perp}\|=\|T \big|_{{H_{0}}^\perp}\|<1$. Now the conclusion follows  from \cite{siju}*{Theorem 2.4}.

(3)$\implies$ (1): Suppose $M_T$ is compact and the norm of $B(H)$ is strongly subdifferentiable at $T$. Now, by \cite{MPSbook}*{Chapter 3, Theorem 3.3.5}, we have $M_T=S_{H_0}$ for some subspace $H_0$ of $H$. Since $M_T$ is compact, it follows that $H_0$ is finite dimensional. As the norm of $B(H)$ is strongly subdifferentiable at $T$, by \cite{siju}*{Theorem 2.4}, we have $\|T|_{{H_0}^{\perp}}\|=\|\,|T|\big|_{{H_0}^{\perp}}\|<1$. Then, by \cite{PSG}*{Theorem 3.1}, we conclude that $T$ satisfies the Bhatia-\v Semrl property.
\end{proof}

In the Hilbert space setting, Theorem~\ref{bshs} shows that if $T \in B(H)$ satisfies the Bhatia-\v{S}emrl property, then $M_T$ is compact and $\spn(M_T)$ is finite dimensional. However, this need not be true in general.
\begin{eg}\label{ceg}
  Let $Y$ be an infinite dimensional reflexive Banach space. Consider the Banach space $X \oplus_{\infty} \mathbb{R}$, which is the direct sum of $X$ and $\mathbb{R}$, equipped with the norm $\|(x, a)\|_{\infty}= \max\{ \|x\|, |a| \}$ for $x \in X$ and $a \in \mathbb{R}$. It is clear that $X$ is a reflexive Banach space and hence, by \cite{KL}*{Corollary 2.2}, we have that every $x^* \in X^*$ satisfies the Bhatia-\v{S}emrl property. Define $f : X \to \mathbb{R}$ as $f(y, a)=a$. Clearly $f \in S_{X^*}$ and $M_f=(B_Y \times \{1\})\cup(B_Y \times \{-1\})$. Since $X$ is reflexive, $f$ satisfies Bhatia-Semrl property. However, $M_f$ is not compact and  $\operatorname{span}(M_f)$ is infinite dimensional.  
\end{eg}
In the next two results, we provide a sufficient condition that relates the Bhatia–\v{S}emrl property of an operator to its essential norm.
\begin{prop}
\label{mtf}
Let $X$ and $Y$ be Banach spaces. Let $T\in B(X, Y)$ be a nonzero operator such that $T$ satisfies the Bhatia-\v Semrl property. If $\spn(M_T)$ is finite dimensional, then $\|T\|_e<\|T\|$.
\end{prop}
 \begin{proof} Assume that $\|T\|_e=\|T\|$. Then $T \perp_B K(X, Y)$. Let $X_0=\spn(M_T)$. Since $X_0$ is a finite dimensional subspace of $X$, it is complemented in $X$.     That is, there exists a  projection $P\in B(X)$
such that $R(P)=X_0$. Hence, $P$ is a finite rank projection.  Now consider the operator $K: X \to Y$ defined as $K=TP$. Clearly, $K\in K(X, Y)$ and hence $T \perp_B K$. 
But for any $x \in M_T$, we have $Kx=TPx=Tx$. Hence, there is no $x \in M_T$ such that $Tx \perp_B Kx$, which contradicts the hypothesis.
 \end{proof}
% \begin{rmk}
%  In view of the Proposition~\ref{mtf}, it is natural to look at the class of operators $T$ on Banach spaces for which $\spn(M_T)$ is finite dimensional. 
% % % In particular, one can ask whether every nonzero operator satisfying the Bhatia-\v Semrl property necessarily has $\operatorname{span}(M_T)$ finite dimensional. 
% % % In the Hilbert space setting, Theorem~\ref{bshs} shows that such operators indeed have $\operatorname{span}(M_T)$ finite dimensional.
%  \end{rmk}
Now let us recall some definitions.
\begin{defn}(\cite{hww})
 A Banach space $X$ is said to have the approximation property (compact approximation property) if, for every $\epsilon>0$ and every compact set $C$ of $X$, there exists a finite rank operator $A \in B(X)$ $($a compact operator $A \in B(X))$ such that $\|x -Ax\|<\epsilon$ for every $x \in C$. 
 % Similarly, $X$ is said to have the compact  approximation property if, for every $\epsilon>0$ and every compact subset $C$ of $X$, there exists an operator $A \in K(X)$ such that $\|x-Ax\|< \epsilon$ for every $x \in C$.
\end{defn}

\begin{prop}\label{es cmt}
    Let $X$ be a Banach space with the compact approximation property, and let $Y$ be any Banach space. If a nonzero operator $T \in B(X, Y)$ satisfies the Bhatia-\v Semrl property and $M_T$ is compact,  then $\|T\|_e<\|T\|$.
\end{prop}
\begin{proof}  Assume that $\|T\|_e=\|T\|$. Then we have $T \perp_B K(X, Y)$. Since $X$ has the compact approximation property, there exists a compact operator $S \in B(X)$ such that $\|x-Sx\| \leq \frac{1}{2}$ for every $x \in M_T$. Let $A=TS$. Then $T\perp_B A$. Thus, for $x \in M_T$, we have $\|Tx-Ax\|=\|Tx-TSx\| \leq \|T\| \|x-Sx\| \leq \frac{\|T\|}{2}<\|T\|=\|Tx\|$.
Thus, $Tx$ is not Birkhoff-James orthogonal $Ax$ for any $x \in M_T$,  which is a contradiction.
\end{proof}

% Motivated by these results, we investigate the relationship among the properties: the Bhatia-\v Semrl property, the condition that the
% essential norm is less than the operator norm, and strong subdifferentiability in $B(\ell_p, \ell_q)$.
We now provide a relation between operators satisfying the Bhatia-\v{S}emrl property and the condition that the essential norm is less than the operator norm.
A Schauder basis $(e_n)$ of a Banach space $X$ is called $1$-unconditional if
$
\left\|\sum_{n=1}^{N} a_n e_n\right\|
\leq
\left\|\sum_{n=1}^{N} b_n e_n\right\|
$
for every $N\in\mathbb{N}$ and all scalars $a_1,\ldots,a_N$ and $b_1,\ldots,b_N$ satisfying $|a_n|\leq |b_n|$ for $1\leq n\leq N$ (see \cite{AK}).
% Recall that a Schauder basis $(e_n)$ in a Banach space $X$ is called $1$-unconditional if, for every $x = \sum_{n=1}^{\infty} a_n e_n$ in $X$, we have $\left\| \sum_{n=1}^{\infty}\epsilon_n a_n e_n \right\| \leq  \|x\|$ for all signs $\epsilon_n= \pm 1$ (see \cite{} for more details).

\begin{prop}\label{bs gen}
    Let $X$ be a Banach space with a $1$-unconditional basis $(e_n)$ and let $Y$ be a Banach space. If a nonzero operator $T \in B(X, Y)$ satisfies the Bhatia-\v{S}emrl property, then $\|T\|_e<\|T\|$.
\end{prop}
\begin{proof}
    Assume that $\|T\|_e=\|T\|$. Then $T \perp_B K(X, Y)$. Define $A : X \to Y$ as  $A(x)= T\left(\sum_{m=1}^{\infty}\frac{1}{m^2}a_me_m \right)$, where $x=\sum_{m=1}^{\infty}a_m e_m$.
 We next show that $A$ is a compact operator. For each $r \in \mathbb{N}$, define $A_r : X \to Y$ as follows, 
$A_r(x)= T\left(\sum_{m=1}^{r}\frac{1}{m^2}a_m e_m \right)$.
For $x=\sum_{m=1}^{\infty}a_m e_m$, we have

\begin{align*}
    \|A(x)-A_r(x)\| &= \left\|T\left(\sum_{m=r+1}^{\infty} \frac{1}{m^2} a_m e_m\right)\right\|\\
    & \leq \|T\| \left\| \sum_{m=r+1}^{\infty} \frac{1}{m^2}a_m e_m \right\|\\
    % & \leq \|T\| \sup_{n \in \mathbb{N}}\left\{\frac{1}{(r+n)^2} \right\}\|x\|\\ 
    & \leq \frac{\|T\|}{(r+1)^2} \|x\|\quad \text{(by \cite{AK}*{Chapter 3, Proposition 3.1.3})}.
\end{align*}

Therefore, $\|A-A_r\| \leq \frac{\|T\|}{(r+1)^2}$ and hence $A_r \to A$ as $r \to \infty$. Thus, $A \in K(X, Y)$ and by assumption, we get $T \perp_B A$.

 We next show that $T$ does not satisfy the Bhatia-\v Semrl property. Let $x=\sum_{m=1}^{\infty}a_m e_m $ be an element in $ M_T$. For $n \in \mathbb{N}$, let $s_n = \sum_{m=1}^{n} a_m e_m$. Then it is clear that $\|s_n-x\| \to 0$. Therefore, there exists $N \in \mathbb{N}$ such that $\|\sum_{m=N+1}^{\infty} a_m e_m\| < \|x\|$. Consider:
\begin{align*}
    \|T(x)-A(x)\|  
   % &= \left\| T\left(\sum_{m=1}^{\infty} x_m e_m\right) - T\left(\sum_{m=1}^{\infty}\frac{1}{m^2}x_m e_m\right) \right\|\\
 &= \left\| T\left(\sum_{m=1}^{\infty} \left(1-\frac{1}{m^2}\right)a_m e_m\right) \right\|\\
    & \leq \|T\|  \left\|\sum_{m=1}^{\infty} \left(1-\frac{1}{m^2}\right) a_m e_m \right\| \\
    % & =\|T\| \left\| \sum_{m=1}^{N}\left(1-\frac{1}{m^2}\right)a_me_m + \sum_{m=N+1}^{\infty}\left(1-\frac{1}{m^2}\right)a_m e_m \right\|\\
    & \leq \|T\| \left\| \sum_{m=1}^{N} \left(1-\frac{1}{N^2}\right)a_m e_m + \sum_{m=N+1}^{\infty} a_me_m \right\| \quad \text{(by \cite{AK}*{Chapter 3, Proposition 3.1.3})}\\
    & = \|T\| \left\| \left(1-\frac{1}{N^2} \right) \sum_{m=1}^{\infty}a_me_m + \frac{1}{N^2} \sum_{m= N+1}^{\infty} a_m e_m \right\| \\
    & \leq  \|T\| \left( \left(1-\frac{1}{N^2}\right)\|x\| + \frac{1}{N^2} \left\| \sum_{m=N+1}^{\infty} a_m e_m \right\| \right)\\
    & < \|T\| \|x\| =\|T(x)\|.
\end{align*}
Thus $T(x)$ is not Birkhoff-James orthogonal to $A(x)$ and hence $T$ does not satisfy the Bhatia-\v Semrl property
\end{proof}

For a net $(x_{\alpha})$ in $X$ and $x \in X$, the weak convergence of $(x_{\alpha})$ to $x$ is denoted by $x_{\alpha} \xrightarrow{w}x$. Recall that a Banach space $X$ has the Kadets-Klee property if, for any sequence $(x_n)$ and $x$ in $X$, the conditions  $x_n \xrightarrow{w} x$ and $\|x_n\| \to \|x\|$ imply $\|x_n - x \| \to 0$.  A sequence $(x_n)$ in $S_X$ is said to be a norming sequence for an operator  $T\in B(X, Y)$ if $\|Tx_n\| \to \|T\|$ and $(x_n)$ is said to be a weakly null norming sequence for $T$ if $x_n \xrightarrow{w} 0$ and $\|Tx_n\| \to \|T\|$. 

 In view of Example \ref{ceg} and Proposition \ref{es cmt}, it is natural to study the compactness of the norm attainment set. We now provide a large class of operators for which the norm attainment set is compact. We use the operator version of property $(M)$, introduced in \cite{KW}. An operator $T$ with $\|T\| \leq 1$ between Banach spaces $X$ and $Y$ is said to have property $(M)$ if for all $x \in X$, $y \in Y$, $\|x \| \leq \|y\|$, and a weakly null sequence $(x_n)$ in $X$, we have $\limsup \|y + Tx_n\| \leq \limsup \|x+x_n\|$.

\begin{prop}\label{mt compact}
    Let $X$ be a reflexive Banach space such that $X$ has the Kadets-Klee property. Let $Y$ be a Banach space. Suppose that $K(X, Y)$ is an $M$-ideal in $B(X, Y)$. If $T \in S_{B(X, Y)}$ is such that $\|T\|_e<1$, then $T$ is norm-attaining and $M_T$ is a compact set.
\end{prop}
\begin{proof}
    Since $\|T\|_e<1$, by \cite{W}*{Lemma 3.1}, it follows  that $M_T \neq \emptyset$. We now show that every norming sequence for $T$ has a subsequence that converges in norm to an element in $M_T$. Consequently, $M_T$ is compact. Let $(x_n)$ be a norming sequence for $T$. Since $X$ is reflexive, there exists a subsequence $(x_{n_{k}})$ of $(x_n)$ such that $x_{n_{k}} \xrightarrow{w}x$. Then $\|x\| = 1$. Suppose, on the contrary, that $\|x\| <1$. Since $K(X, Y)$ is an $M$-ideal in $B(X, Y)$, from the proof of \cite{KW}*{Theorem 6.3} and \cite{DW}*{Proposition 2.3 and from the proof of $1 \implies 2$ of Theorem 3.1}, we conclude that every contraction in $B(X, Y)$ has the property $(M)$. Now let $r>0$ be such that $\|T\|_e<r<1$. Then there exists a compact operator $K \in K(X, Y)$ such that $\|T-K\|<r$. Since $\frac{T-K}{r}$ is a contraction and $(r (x_{n_{k}}-x))$ is a weakly null sequence, property $(M)$ yields
    
    $$\limsup \left\|Tx + \left(\frac{T-K}{r}\right)(r(x_{n_{k}}-x))\right\| \leq \limsup\|x+r(x_{n_k}-x)\|.$$
    Therefore, we  have $$\limsup \|Tx_{n_{k}}-K(x_{n_{k}}-x)\| \leq \limsup\|(1-r)x+rx_{n_{k}}\|.$$ 
    Observe that $$ \Big|\|Tx_{n_{k}}\| - \|Tx_{n_{k}} - K(x_{n_{k}}-x)\| \Big|\leq \|K(x_{n_k}-x)\| \to 0.$$  Hence, we get $$1= \limsup \|Tx_{n_{k}}-K(x_{n_{k}}-x)\| \leq \limsup\|(1-r)x+rx_{n_{k}}\|\leq (1-r)\|x\|+r<1,$$ which is a contradiction. Therefore $\|x\|=1$. Since $X$ has the Kadets-Klee property, we get $\|x_{n_{k}}-x\| \to 0$ and $x \in M_T$. 
    % As a consequence, we conclude that $M_T$ is compact.
    % Let $(z_n)$ be a sequence in $M_T$. Then $\|Tz_n\|=1$. Hence, by our earlier argument, there exists $z \in M_T$ and a subsequence $(z_{n_k})$ of $(z_n)$ such that $\|z_{n_{k}} -z\| \to 0$. Therefore, $M_T$ is compact.
\end{proof}

% Recall that, for $T\in B(X, Y)$, a sequence $(x_n)$ in $S_X$ is said to be a norming sequence if $\|Tx_n\| \to \|T\|$, and $(x_n)$ is said to be a weakly null norming sequence for $T$ if $x_n \xrightarrow{w} 0$ and $\|Tx_n\| \to \|T\|$.
% \begin{lem}\label{ess weakly null}
%     Let $X$ and $Y$ be Banach spaces. If $T \in B(X, Y)$ has a weakly null norming sequence, then $\|T\|_e=\|T\|$.
% \end{lem}
% \begin{proof}
%     Let $(x_n)$ be a sequence in $S_X$ such that $x_n \xrightarrow{w}0$ and $\|Tx_n\| \to \|T\|$. Now let $L \in K(X, Y)$. Then $\|T+L\|\geq \|Tx_n+Lx_n\|\geq \|Tx_n\|-\|Lx_n\|$. Since $Lx_n \to 0$ in the norm topology, we have $\|T\|\leq \|T+L\|$ and hence $\|T\| \leq \|T\|_e$. Therefore $\|T\|_e=\|T\|$.
% \end{proof}

The following result is a  consequence of the proof of \cite{hww}*{Chapter \Rmnum{6}, Proposition 4.7}. 
\begin{prop}\label{weakly null chara}
    Let $X$ and $Y$ be reflexive Banach spaces such that $K(X,Y)$ is an M-ideal in $B(X,Y)$.  Let $T \in B(X,Y)$ be a nonzero operator. Then $\|T\|_e=\|T\|$ if and only if either $T$ or $T^*$ has a weakly null norming sequence.
\end{prop}
\begin{proof}
Let $(x_n)$ be a weakly null norming sequence for $T$. Then for $L \in K(X, Y)$, we have  $\|T+L\|\geq \|Tx_n+Lx_n\|\geq \|Tx_n\|-\|Lx_n\|$. Since $Lx_n \to 0$ in the norm topology, we have $\|T\|\leq \|T+L\|$ and hence $\|T\| \leq \|T\|_e$. Therefore $\|T\|_e=\|T\|$.
 Now suppose $T^*$ has a weakly null norming sequence.  Then, $\|T^*\|_e=\|T^*\|=\|T\|$. Since $X$ and $Y$ are reflexive, it follows from \cite{M}*{Chapter 3, Theorem 3.1.11} that $\|T^*\|_e=\|T\|_e$. Therefore, $\|T\|_e=\|T\|$.

Conversely, suppose $\|T\|_e=\|T\|$, then it follows from the  proof of \cite{hww}*{Chapter \Rmnum{6}, Proposition 4.7} that  either $\|T\|_e=\limsup\|Tx_{\alpha}\|$ for some net $(x_{\alpha})$ in $S_X$ with $x_{\alpha} \xrightarrow{w} 0$ or $\|T\|_e=\limsup\|T^*y_{\alpha}^*\|$ for some net $(y_{\alpha}^*)$ in $S_{Y^*}$ with $y_{\alpha}^* \xrightarrow{w}0$. Now, suppose $\|T\|_e=\limsup \|Tx_{\alpha}\|$ for some net $(x_\alpha)_{\alpha \in I}$ in $S_X$ such that $x_{\alpha} \xrightarrow{w}0$. Consider the Banach space $X \oplus_{\infty} \mathbb{R}$, which is the direct sum of $X$ and $\mathbb{R}$, equipped with the norm $\|(x, a)\|_{\infty}= \max\{ \|x\|, |a| \}$ for $x \in X$ and $a \in \mathbb{R}$. Since $X$ is reflexive, it is easy to see that $B_X \times [0, \|T\|]$ is weakly compact in $X \oplus_{\infty}\mathbb{R}$. Consider the set $A =\{(x_{\alpha}, \|Tx_{\alpha}\|) : \alpha \in I\} \subseteq B_X \times [0,\|T\|]$.  Then the weak closure of $A$, denoted by $\overline{A}^w$, is weakly compact in $X \oplus_{\infty}\mathbb{R}$ and $(0, \|T\|) \in \overline{A}^w$. Therefore, by \cite{M}*{Chapter 2, Lemma 2.8.5}, we can find a sequence $((x_{\alpha_{n}}, \|Tx_{\alpha_{n}}\|))$ in $A$ such that $(x_{\alpha_{n}}, \|Tx_{\alpha_{n}}\|) \xrightarrow{w}(0, \|T\|)$. Hence $x_{\alpha_{n}} \xrightarrow{w}0$ and $\|Tx_{\alpha_{n}}\| \to \|T\|$. Similarly, we can show that $T^*$ has a weakly null norming sequence if $\|T\|_e=\limsup\|T^*y_{\alpha}^*\|$ for some net $(y_{\alpha}^*)$ in $S_{Y^*}$ with $y_{\alpha}^* \xrightarrow{w}0$.  
\end{proof}
\begin{rmk}\label{ssdrmk}
In \cite{siju}*{Remark 2.12}, the points of strong subdifferentiability of the norm of $B(X,Y)$ were characterized under the assumptions that $X$ is uniformly convex, $Y$ is uniformly smooth, $K(X,Y)$ is an $M$-ideal in $B(X,Y)$, $X^*$ has the Radon–Nikodým property, and either $X^*$ or $Y^*$ has the approximation property. We observe that the last two assumptions can be omitted. Indeed, the approximation-property assumption is used in the proof of \cite{siju}*{Theorem 2.11} only to identify $K(X, Y)^*$ with the projective tensor product $X \hat{\otimes}_{\pi}Y^*$
 However, the argument in the proof of \cite{C}*{Theorem 3.3} yields this identification whenever $X$ and $Y$ are reflexive and $K(X,Y)$ is an $M$-ideal in $B(X,Y)$. Moreover, the reflexivity of $X$ implies that $X^*$ has the Radon–Nikod\'ym property. Consequently, the characterization in \cite{siju}*{Remark 2.12} remains valid under the sole assumptions that $X$ is uniformly convex, $Y$ is uniformly smooth, and $K(X,Y)$ is an $M$-ideal in $B(X,Y)$.
\end{rmk}

    In \cite{siju}*{Theorem 3.1}, it was shown that if $T \in B(\ell^p, \ell^q)$ with $1<p, q<\infty$ satisfies $\|T\|_e<\|T\|$, then the norm of $B(\ell^p, \ell^q)$ is strongly subdifferentiable at $T$. However, the converse need not be true. Indeed, for any Banach space $X$, by \cite{MJAP}*{Proposition 4.5} and \cite{FP}*{Theorem 1.2}, the norm of $B(X)$ is strongly subdifferentiable at the identity operator $I$, but $\|I\|_e=\|I\|=1$. The next result shows that the converse holds for certain classes of Banach spaces $X$ and $Y$ when $M_T$ is compact.
\begin{prop}\label{ess char}
    Let $X$ be a uniformly convex Banach space and  $Y$ be a uniformly smooth Banach space such  that $K(X, Y)$ is an $M$-ideal in $B(X, Y)$. Let $T \in B(X, Y)$ be a nonzero operator. Then $\|T\|_e<\|T\|$ if and only if $M_T$ is compact and the norm of $B(X, Y)$ is strongly subdifferentiable at $T$.
\end{prop}
\begin{proof}
    Let $T \in B(X, Y)$ be such that $\|T\|_e<\|T\|$. Then, by Proposition \ref{mt compact}, $M_T$ is a nonempty compact subset of $X$. Now let $(x_n)$ be a norming sequence for $T$. Then, from the proof of Proposition \ref{mt compact}, we obtain a subsequence $(x_{n_k})$ of $(x_n)$ and an element $x \in M_T$ such that $\|x_{n_k}-x\| \to 0$ and hence, $d(x_{n_{k}}, M_T) \to 0$. Now, by Remark \ref{ssdrmk}, we conclude that the norm of $B(X, Y)$ is strongly subdifferentiable at $T$.
    
      Conversely, let $T \in B(X, Y)$ be such that the norm of $B(X, Y)$ is strongly subdifferentiable at $T$ and $M_T$ is compact. Then, by Remark \ref{ssdrmk}, we have that $M_T \neq \emptyset$. Assume $\|T\|_e=\|T\|$. Hence, by Proposition \ref{weakly null chara}, it follows that  either $T$ or $T^*$ has a weakly null norming sequence. 
\begin{case}
  Suppose $T$ has a weakly null norming sequence, say $(x_n)$.  Since the norm of $B(X, Y)$ is strongly subdifferentiable at $T$, by Remark \ref{ssdrmk},  we can conclude that $(x_n)$ has a subsequence $(x_{n_{k}})$ such that $d(x_{n_{k}}, M_T) \to 0$. 
  Thus, we can find a sequence $(z_{n_{k}})$ in $M_T$ such that $\|x_{n_{k}}-z_{n_{k}}\| \to 0$. Since $M_T$ is compact, $(z_{n_{k}})$ has a subsequence, which we denote again by $(z_{n_{k}})$, such that $\|z_{n_{k}} - z\| \to 0$ for some $z \in M_T$. Now consider $\|x_{n_{k}}-z\| \leq \|x_{n_{k}} - z_{n_{k}}\|+\|z_{n_{k}} - z\| \to 0$. 
  Since $x_{n_{k}} \xrightarrow{w} 0$, we get $z=0$, which is a contradiction. 
\end{case}

\begin{case}
    Suppose $T^*$ has a weakly null norming sequence. Let $(y^*_n)$ be a sequence in $S_{Y^*}$ such that $y^*_n \xrightarrow{w}0$ and $\|T^*y^*_n\| \to \|T\|$. Since $X$ is reflexive, each $T^*y^*_n$ attains its norm. Let $x_n \in S_X$ be such that $T^*y^*_n(x_n)=\|T^*y^*_n\|$. Therefore, $\|Tx_n\| \to \|T\|$. Since the norm of $B(X, Y)$ is strongly subdifferentiable at $T$, by Remark \ref{ssdrmk}, we conclude that $(x_n)$ has a subsequence $(x_{n_{k}})$ such that $d(x_{n_{k}}, M_T) \to 0$. Thus, we can find a sequence $(z_{n_{k}})$ in $M_T$ such that $\|x_{n_{k}} - z_{n_k}\| \to 0$. Since $M_T$ is compact, $(z_{n_k})$ has a subsequence, which we denote again by $(z_{n_{k}})$, such that $\|z_{n_{k}} - z\| \to 0$ for some $z \in M_T$. Now, consider $\|x_{n_{k}}-z\| \leq \|x_{n_{k}} - z_{n_{k}}\|+\|z_{n_{k}} - z\| \to 0$. Hence, $\|Tx_{n_{k}} - Tz\| \to 0$. Since $y^*_{n_{k}} \xrightarrow{w}0$, we have $y^*_{n_{k}}(Tx_{n_{k}}) \to 0$, which is a contradiction to $y^*_{n_{k}}(Tx_{n_{k}}) = \|T^*y^*_{n_{k}}\| \to \|T\|$. \qedhere
\end{case}
 \begin{rmk}
 It follows from \cites{hww, KW} that the conclusion of Proposition \ref{ess char} holds for the following pairs $(X, Y)$ of Banach spaces: $(\ell^p, \ell^q)$ where $1<p, q< \infty$, $(\ell^p, L^p[0, 1])$ where $2\le p < \infty$ and $(\ell^p, c_p)$ where $2\leq p< \infty$ and $c_p$ denote the Schatten class of all compact operators on $\ell^2$ with the $p$-summable singular values. 
 % See \cites{hww, KW}, where it is shown that if $(X, Y)$ is one of the above mentioned pairs of Banach spaces, then $K(X, Y)$ is an $M$-ideal in $B(X, Y)$.
    % For $1<p, q<\infty$, it follows from \cite{hww}*{Chapter \Rmnum{6}, Example 4.1} that $K(\ell^p, \ell^q)$ is an $M$-ideal in $B(\ell^p, \ell^q)$. Moreover, for $2\leq p < \infty$, it follows from \cite{KW}*{Proposition 6.5} that $K(\ell^p, L^p[0.1])$ is an $M$-ideal in $B(\ell^p, L^p[0,1])$. Furthermore, let $c_p$ denote the Schatten class of all compact operators on $\ell^2$ with the $p$-summable singular values. Hence, for $2\leq p<\infty$, it follows from \cite{KW}*{Proposition 6.6} that $K(\ell^p, c_p)$ is an $M$-ideal in $B(\ell^p, c_p)$. These are some of the classical Banach spaces for which Proposition \ref{ess char} holds true.
 \end{rmk}
\end{proof}
    The next result follows immediately from Proposition \ref{bs gen} and Proposition \ref{ess char}.
\begin{cor}\label{bs ess}
Let $Y$ be a uniformly smooth Banach space. For $1<p<\infty$, suppose that $K(\ell^p, Y)$ is an $M$-ideal in $B(\ell^p, Y)$. Let $T \in B(\ell^p, Y)$ be such that $T\neq 0$. Consider the following statements:
 \begin{enumerate}
     \item $T$ satisfies the Bhatia-\v Semrl property.
     \item $\|T\|_e<\|T\|$.
     \item $M_T$ is compact, and the norm of $B(\ell^p, Y)$ is strongly subdifferentiable at $T$.
 \end{enumerate}
 Then the following implications hold: $(1) \implies (2)$ and $(2) \iff (3)$.
\end{cor}

As an application of Corollary \ref{bs ess}, we give a class of operators that satisfy the Bhatia-\v Semrl property. For $1<p<\infty$, let $T \in B(\ell^p)$ be a diagonal operator defined as $T(x_n)=(k_n x_n)$, where $(k_n)$ is a real valued bounded sequence. Then it is well known that $\|T\|_e = \limsup |k_n|$.
\begin{prop}
    \label{diag norm}
 For $1 <p<\infty $, let $T \in  S_{B(\ell^p)}$ be the diagonal operator defined as $T(x_n)= (k_nx_n)$, where $(k_n)$ is a real valued bounded sequence. Then $T$ satisfies the Bhatia-\v Semrl property if and only if $\|T\|_e < 1$. Consequently, every compact diagonal operator on $\ell^p$ satisfies the Bhatia-\v Semrl property.
 \end{prop}
 \begin{proof}
  %Suppose $M_T \neq \emptyset$ and let $x\in M_T$. Let $J: \ell_p \to \ell_p^*$ be the duality map. Then $J(x)= (|x_n|^{p-1}\sgn(x_n))$. We observe that $J(Tx)=\frac{1}{\|T\|^{p-2}}(|k_n x_n|^{p-1}\sgn(k_nx_n))$. Thus 
  %\begin{align*}
  %T^*J(Tx) =\frac{1}{\|T\|^{p-2}}(|k_nx_n|^{p-1}k_n \sgn(k_nx_n))
  %&=\frac{1}{\|T\|^{p-2}}(|k_nx_n|^{p-1}k_n \sgn(k_n)\sgn(x_n))\\
 %& =\frac{1}{\|T\|^{p-2}}(|k_n|^p|x_n|^{p-1}\sgn(x_n)).
  %\end{align*}
   Suppose $T$ satisfies the Bhatia-\v Semrl property. Then, by Corollary \ref{bs ess}, we get $\|T\|_e < 1$. 

 Conversely, suppose $\|T\|_e<1$. Then, by Proposition \ref{mt compact}, $T$ attains its norm.  We first show that the set $\{|k_n| : n \in \mathbb{N}\}$ attains its supremum value. Let $x \in M_T$. Then we have $\|Tx\|^p=\sum_{n=1}^{\infty}|k_nx_n|^p=1=\sup\{|k_n| : n \in \mathbb{N}\}$. Therefore, whenever $x_n \neq 0$, we have $|k_n|=1$. Hence, the set $\{|k_n| : n \in \mathbb{N}\}$ attains its supremum value. If the set $\{|k_n| : n \in \mathbb{N}\}$ attains its supremum at infinitely many $n_m$, then $(e_{n_m})$ would be a weakly null norming sequence for $T$.  Hence, from the proof of Proposition \ref{weakly null chara}, we get $\|T\|_e=1$, a contradiction.  So let $n_1, n_2, \dots, n_m$ be all the indices at which $\{|k_n| : n \in \mathbb{N}\}$ attains its supremum value. Let $X_0 =\spn \{e_{n_{i}} : 1\leq i \leq m \}$. Then $M_T=S_{X_0}$. Indeed, let $x \in M_T$. Then we have $\|Tx\|^p=\sum_{n=1}^{\infty}|k_nx_n|^p=1$. Therefore, $\sum_{n=1}^{\infty}(1-|k_n|^p)|x_n|^p=0$. Now note that for every $i \in \mathbb{N} \setminus\{n_1, n_2, \dots, n_m\}$, we have $1-|k_i|^p>0$. Hence, $x_i=0$ for every $i \in \mathbb{N} \setminus\{n_1, n_2, \dots, n_m\}$. Therefore, $x \in S_{X_0}$. Now let $x \in S_{X_0}$.
Then $x=\sum_{i=1}^{m} a_{n_{i}} e_{n_{i}}$ for some scalars $a_{n_{i}}$. 
 %Let $x, y \in M_T$ and let $a, b \in \mathbb{R}$, 
Therefore
\begin{align*}
   \left \|T\left(\sum_{i=1}^{m} a_{n_{i}} e_{n_{i}}\right) \right \|_p^p  = \left \| \sum_{i=1}^{m} k_{n_{i}} a_{n_{i}} e_{n_{i}} \right \|_p^p
    =\sum_{i=1}^{m} |k_{n_{i}}a_{n_{i}} |^p 
    = \sum_{i=1}^{m}|a_{n_i}|^p
    =1.
\end{align*}
 Hence $M_T=S_{X_0}$. Then, by \cite{MPRS}*{Theorem 3.2}, we can see that $T$ satisfies the Bhatia-\v Semrl property.
 \end{proof}

 \section{Bhatia-Semrl property of norm one projections.}
 \label{bjp}
In this section, we discuss the Bhatia-\v{S}emrl property of the norm one projections on Banach spaces. In general, a norm one projection on a Banach space $X$ need not satisfy the Bhatia-\v{S}emrl property. For example, it follows from \cite{CK}*{Proposition 3.7} and \cite{KL}*{Proposition 3.6} that norm one projections on $c_0$ and $L^1[0, 1]$ do not satisfy the Bhatia-\v Semrl property, respectively. However, we identify conditions on a Banach space $X$ under which every finite rank norm one projection on $X$ satisfies the Bhatia-\v{S}emrl property. For a projection $P$ on a Banach space $X$,  the range and kernel of $P$ are denoted by  $R(P)$ and $\ker(P)$, respectively.

 \begin{prop}\label{finte rank bs}
     Let $X$ be a reflexive, strictly convex Banach space with the Kadets-Klee property. Then any finite rank norm one projection on $X$ satisfies the Bhatia-\v Semrl property.
\end{prop}
\begin{proof}
Let $P$ be a norm one projection on $X$. We first show that $M_P=S_{R(P)}$. Since $\|P\|=1$, it is clear that $S_{R(P)} \subseteq M_P$. Let $x \in M_P$ and let $y=Px$. Hence, $P\left(\frac{x+y}{2}\right)=y$. If $x \neq y$, then by the strict convexity of $X$, we get that $\left\|\frac{x+y}{2}\right\|<1$. Therefore, we have $\|P\| \geq \frac{\left\|P\left(\frac{x+y}{2}\right)\right\|}{\left\| \frac{x+y}{2}\right\|}>1$, which is not possible. Hence, $x=y$ and $M_P\subseteq R(P)$. Since $R(P)$ is finite dimensional, we have that $M_P$ is compact. Next, let $C$ be a closed subset of $S_X$ such that $d(C, M_P)>0$.  Assume that $\sup \{\|Px\| : x \in C \}=1$. Then there exists a sequence $(x_n)$ in $C$ such that $\|Px_n\| \to 1$. Since $X$ is reflexive, there exists a subsequence $(x_{n_{k}})$ of $(x_n)$ such that $x_{n_{k}} \xrightarrow{w} x_0$, for some $x_0 \in B_X$. Since $P$ is compact, we have $Px_{n_{k}} \to Px_0$ in norm. Hence, $\|Px_0\|=1$.
 Thus, we have $x_0 \in M_P \subset S_X$. Since $X$ has the Kadets-Klee property, we have $x_{n_{k}} \to x_0$ in norm. Thus, $d(C, M_P) \leq \|x_{n_{k}} -x_0 \| \to 0$, which is a  contradiction to the choice of $C$. Therefore, $\sup \{ \|Px\|  : x \in C \} <1$. Hence, by \cite{PSG}*{Theorem 2.2}, 
    it follows that $P$ satisfies the Bhatia-\v Semrl property.
\end{proof}
 % The following remark is useful in characterizing the Bhatia-\v{S}emrl property of norm one projections on $\ell^p$.

\begin{rmk}\label{finite dim proj}
   One can easily observe that for a Banach space $X$, a norm one projection $P \in B(X)$ has finite rank if and only if  $\|P\|_e<1$. For, suppose $\|P\|_e<1$ and rank of $P$ is  infinite. Then there exists  a $\delta>0$ and an operator $A \in K(X)$ such that $\|P\|_e\leq \|P-A\|<\delta<1$. Now define $L : R(P) \to R(P)$ as $L= PA$. Then $L$ is a compact operator and for $x \in S_{R(P)}$, we have $\|x-Lx\|=\|P^2x-PAx\|\leq \|P-A\|\|x\|<\delta$. Therefore $\|I-L\|\leq \delta <1$ and hence $L$ is an invertible compact operator on $R(P)$, which is a contradiction.
\end{rmk}
 The next result follows immediately from Proposition \ref{finte rank bs}, Remark \ref{finite dim proj} and Corollary \ref{bs ess}.
\begin{cor}
For $1<p<\infty$, let $P \in B(\ell^p)$ be a norm one projection. Then $P$ satisfies the Bhatia-\v Semrl property if and only if $P$ is finite rank.
\end{cor}

We now relax the conditions on the Banach space $X$ and prove a result similar to Proposition \ref{finte rank bs} for $L^p$-Projections $(1\leq p < \infty)$.

 \begin{prop}\label{LP BS}
    Let $X$ be a Banach space and let $P \in B(X)$ be a finite rank $L^p$-projection. Then $M_P=S_{R(P)}$ and $P$ satisfies the Bhatia-\v Semrl property.
 \end{prop}
\begin{proof}  
Clearly $S_{R(P)} \subseteq M_P$. Now let $x \in M_P$. Then $1=\|x\|^p=\|Px\|^p+\|x-Px\|^p=1+\|x-Px\|^p$. Thus $x=Px \in S_{R(P)}$. Hence $M_P=S_{R(P)}$. Note that since $P$ is of finite rank, $S_{R(P)}$ is compact.  Let $C$ be a closed subset of $S_X$  such that $d(C,  S_{R(P)})>0$, where $d(C, S_{R(P)})=\inf \{\|y-z\| : y \in C \mbox{ and } z \in S_{R(P)}\}$. Assume that $\sup \{\|Px\| : x \in C \}=\|P\|=1$,  then there exists a sequence $(x_n)$ in $C$
 such that $\|Px_n\| \to \|P\|=1$. Since $P$ is an $L^p$-projection, we get  $\|x_n- Px_n\| \to 0$. Note that $\|Px_n\| \leq 1$. Since $R(P)$ is finite dimensional, there exists a subsequence $(Px_{n_k})$ of $(Px_n)$,  such that $Px_{n_{k}}\to y$. Since $R(P)$ is closed and  $\|Px_{n_k}\| \to 1$, we get $y \in S_{R(P)}$.  However, we note that $d(C, S_{R(P)}) \leq \|x_{n_k}-y\| \leq \|y-Px_{n_k}\| + \|Px_{n_k} -x_{n_k}\| \to 0$, which is a contradiction to the choice of $C$. Therefore $\sup\{\|Px\| : x \in C \}<1$.
Hence, by \cite{PSG}*{Theorem 2.2}, it follows that $P$ satisfies the Bhatia-\v Semrl property.
\end{proof}

\begin{rmk}
  A result similar to Proposition \ref{LP BS} does not hold in general for $M$-projections. For,  finite rank $M$-projections on $c_0$ do not satisfy the Bhatia-\v Semrl property (see \cite{CK}*{Proposition 3.7}).
\end{rmk}

 We next prove the converse of Proposition \ref{LP BS} for the space $L^p(\mu)$. We first make the following remark.

\begin{rmk}\label{Infinite LP Proj}
    It follows from \cite{BDE} that if $(\Omega, \mathcal{F},\mu)$ is a $\sigma$-finite measure space, then for any $L^p$-projection $P$ on $L^p(\mu)$ $(1\leq p< \infty \; \text{with} \; p\neq 2)$ , there exists a measurable set $D $ such that $P(f)=\chi_D f$ for every $f\in L^p(\mu)$. Moreover, if $P$ has infinite rank, then  there exists a mutually disjoint sequence $(D_n)$ of measurable subsets of $D$ such that $0<\mu(D_n)<\infty$ for all $n$ and $D=( \cup_{n=1}^{\infty}D_n)\cup N$ for some measure zero set $N$. Indeed, let $D=A\cup B$, where $A$ is the union of all atoms of $\mu$ contained in $D$ and $B=D\setminus A$. Since $\mu$ is $\sigma$-finite, the required conclusion follows  if $D$ has infinitely many atoms. So suppose $D$ contains only finitely many atoms. Since $P$ is a projection of  infinite rank and since every $f \in L^p(\mu)$ is almost everywhere constant on atoms, we get $\mu(B)>0$. Since $\mu$ is $\sigma$-finite, there exists  $C \subset B$ such that $0<\mu(C)<\infty$. Since $C$ is not an atom, there exists a $C_1\subset C$ such that $0<\mu(C_1)<\infty$ and $0<\mu(C\setminus C_1)<\infty$. Repeat the argument with $C\setminus C_1$ and then proceeding  similarly, we get a mutually disjoint sequence $(C_n)$ of subsets of $C$ such that $0<\mu(C_n)<\infty$. Then $D=A\cup (\cup_{n=1}^{\infty}C_n)\cup (B\setminus (\cup_{n=1}^{\infty}C_n))$. Now the required conclusion follows, as $B\setminus (\cup_{n=1}^{\infty}C_n)$ can be written as the countable disjoint union of measurable sets of finite measure.

    \end{rmk}

    \begin{prop}\label{lp proj bs}
 Let $(\Omega, \mathcal{F}, \mu)$ be a positive $\sigma$-finite measure space, and let $P : L^p(\mu) \to L^p(\mu)$ be an infinite rank $L^p$-projection. Then $P$ does not satisfy the Bhatia-\v Semrl property.
\end{prop}
\begin{proof}
Let $P$ be an infinite rank $L^p$-projection. For $p=2$, the result follows from Theorem \ref{bshs}. We therefore assume that $p\neq 2$. Then, there exists $D \in \mathcal{F}$ such that $P(f)= \chi_D f$ for $f \in L^p(\mu)$. From Remark \ref{Infinite LP Proj}, we can find a mutually disjoint sequence $(D_n)$ of measurable subsets of $D$ such that $0<\mu(D_n)<\infty$ for every $n$ and $D= (\cup_{n=1}^{\infty}D_n) \cup N$ for some measure zero set $N$.

Now for $f \in R(P)$, we have
   \begin{align}
   \label{frange}
    \sum_{n=1}^{\infty} \int_{\Omega} | \chi_{D_{n}}f |^p d\mu  = \sum_{n=1}^{\infty} \int_{D_n}|f|^p d\mu 
     = \int_{\cup_{n=1}^{\infty}D_n} |f|^p d\mu 
     = \int_{D} |f|^p d\mu 
     = \int_{\Omega} |f|^p d\mu 
    = \|f\|^p_p.
   \end{align}
   
   Define $A : L^p(\mu) \to L^p(\mu)$ as $A(f) = \sum_{k=1}^{\infty}\frac{1}{2^k} \chi_{D_{k}}f$. Then, for $f \in L^p(\mu)$, we have
   \begin{align}
        \|A(f)\|_p  = \left\| \sum_{k=1}^{\infty} \frac{1}{2^k} \chi_{D_{k}} f\right\|_p 
    \leq \sum_{k=1}^{\infty} \frac{1}{2^k} \| \chi_{D_{k}} f\|_p 
    \leq \left(\sum_{k=1}^{\infty}\frac{1}{2^k}\right)\|f\|_p.
   \end{align}
Thus, $A$ is a bounded linear operator on $L^p(\mu) $. 
 
We next show that $P \perp_B A$. For $\lambda \in \mathbb{R}$, we have 
\begin{align*}
    \|P+ \lambda A\|  \geq \left\|P\left(\frac{\chi_{D_{m}}}{\| \chi_{D_{m}}\|_p}\right) + \lambda A\left(\frac{\chi_{D_{m}}}{\| \chi_{D_{m}}\|_p}\right) \right\|_p 
     = \left\| \frac{\chi_{D_{m}}}{\| \chi_{D_{m}}\|_p} +\lambda \left(\frac{\chi_{D_{m}}}{2^m \| \chi_{D_{m}}\|_p}\right)\right\|_p 
     =\left|1+\frac{\lambda}{2^m}\right| .
\end{align*}
Thus we have $\|P+ \lambda A\| \geq \left|1+\frac{\lambda}{2^m}\right| \to 1=\|P\|$ and hence $P \perp_B A$.
Now observe that  for $f\in L^p(\mu)$ and  for a fixed $n$, we have $\chi_{D_{n}}(f-\sum_{k=1}^{\infty}\frac{1}{2^k} \chi_{D_{k}}f)= (1-\frac{1}{2^n})\chi_{D_{n}}f$ almost everywhere on $\Omega$. Note that $M_P=S_{R(P)}$.  Let $f \in M_P$, then $P(f)-A(f)= f-\sum_{k=1}^{\infty}\frac{1}{2^k}\chi_{D_{k}}f \in R(P)$. Hence, by (\ref{frange}), we have 
\begin{align*}
  \|P(f) -A(f) \|_p^p  = \left\| f - \sum_{k=1}^{\infty} \frac{1}{2^k}\chi_{D_{k}}f \right\|_p^p 
  & =  \sum_{n=1}^{\infty} \int_{\Omega} \left| \chi_{D_{n}} \left(f- \sum_{k=1}^{\infty} \frac{1}{2^k}\chi_{D_{k}}f\right) \right|^p d\mu \\
  & = \sum_{n=1}^{\infty} \int_{\Omega} \left| \left(1-\frac{1}{2^n}\right) \chi_{D_{n}} f \right|^p d\mu \\
  & = \sum_{n=1}^{\infty} \left(1-\frac{1}{2^n}\right)^p \int_{D_n}|f|^p d\mu \\
  & < \sum_{n=1}^{\infty} \int_{D_n} |f|^p d\mu \\
  & \leq\|f\|^p_p
   = \|P(f)\|_p^p.
\end{align*}
Therefore, $P(f)$ is not Birkhoff-James orthogonal to $A(f)$.
\end{proof}

It follows from \cite{PSG}*{Theorem 2.1} that isometries on finite dimensional Banach spaces satisfy the Bhatia-\v Semrl property. In contrast, we end this section by proving that if $L^p(\mu)$ is infinite dimensional, then no isometry between $L^p(\mu)$ satisfies the Bhatia-\v{S}emrl property. Since the proof of this is similar to Proposition \ref{lp proj bs}, we just provide a sketch of the proof. 
% In contrast, if $(\Omega, \mathcal{F}, \mu)$ is a $\sigma$-finite positive measure space such that $L^p(\mu)$  is infinite dimensional, then we show that isometries between $L^p(\mu)$ do not satisfy the Bhatia-\v Semrl property.
\begin{prop}
Let $(\Omega, \mathcal{F}, \mu)$ be a positive $\sigma$-finite measure space such that $L^p(\mu)$ $(1\leq p<\infty)$ is infinite dimensional. Then no isometry from $L^p(\mu)$ to $L^p(\mu)$ satisfies the Bhatia-\v Semrl property.  
\end{prop}
\begin{proof}
For $p=2$, the result follows from Theorem \ref{bshs}. We therefore assume that $p\neq 2$.
    We first observe that $\Omega$ cannot be a disjoint union of finitely many atoms as that would make $L^p(\mu)$ finite dimensional. Thus, using an argument similar to Remark \ref{Infinite LP Proj}, we have $\Omega =(\cup_{k=1}^{\infty} D_k) \cup N$,  where $(D_n)$ is a sequence of mutually disjoint measurable sets such that $0 < \mu(D_k)< \infty$ and $N$ has measure zero.
    Let $T \in B(L^p(\mu))$ be an isometry. 
    
    Define $A : L^p(\mu) \to L^p(\mu)$ as $A(f)= T(\sum_{k=1}^{\infty}\frac{1}{2^k}\chi_{D_{k}}f)$. Then proceeding as in the proof of Proposition~\ref{lp proj bs}, we can see that $A\in B(L^p(\mu))$, $T \perp_B A$ and $T(f)$ is not Birkhoff-James orthogonal to $A(f)$ for any $f\in S_{L^p(\mu)}$. Therefore, $T$ does not satisfy the Bhatia-\v Semrl property.
    \end{proof}

\section*{Acknowledgements.}
The authors would like to thank  Sreejith Siju for the Remark \ref{finite dim proj}.
\section*{Disclosure statement.}
No potential conflict of interest was reported by the authors.
\section*{Funding.}
The research of the first named author is supported by KSCSTE, Kerala through the KSYSA-Research Grant and the research of the second named author is supported by the Prime Minister's Research Fellowship (PMRF ID - 3103662).
% \section*{Declaration of generative AI use.}
% The authors report that generative AI was not used in their research or in the preparation of this manuscript.

\bibliography{bibliography.bib}

@book {MPSbook,
    AUTHOR = {Mal, Arpita and Paul, Kallol and Sain, Debmalya},
     TITLE = {Birkhoff-{J}ames orthogonality and geometry of operator
              spaces},
    SERIES = {Infosys Science Foundation Series},
      NOTE = {},
 PUBLISHER = {Springer, Singapore},
      YEAR = {2024},
     PAGES = {xiii+251},
      ISBN = {978-981-99-7110-7; 978-981-99-7111-4},
   MRCLASS = {46-02 (46A32 46B20 46B28 46C50 47A58 47L05 47L25)},
}

@article {DW,
    AUTHOR = {Werner, Dirk},
     TITLE = {{$M$}-ideals and the ``basic inequality''},
   JOURNAL = {J. Approx. Theory},
  FJOURNAL = {Journal of Approximation Theory},
    VOLUME = {76},
      YEAR = {1994},
    NUMBER = {1},
     PAGES = {21--30},
      ISSN = {0021-9045,1096-0430},
   MRCLASS = {47D15 (46B28 47A58 47B07 47D25)},
   }

@article {W,
    AUTHOR = {W\'ojcik, Pawe\l},
     TITLE = {Birkhoff orthogonality in classical {$M$}-ideals},
   JOURNAL = {J. Aust. Math. Soc.},
  FJOURNAL = {Journal of the Australian Mathematical Society},
    VOLUME = {103},
      YEAR = {2017},
    NUMBER = {2},
     PAGES = {279--288},
      ISSN = {1446-7887,1446-8107},
   MRCLASS = {46B20 (46C50 47L05)},
   }

@article {KW,
    AUTHOR = {Kalton, Nigel J. and Werner, Dirk},
     TITLE = {Property {$(M)$}, {$M$}-ideals, and almost isometric structure
              of {B}anach spaces},
   JOURNAL = {J. Reine Angew. Math.},
  FJOURNAL = {Journal f\"ur die Reine und Angewandte Mathematik. [Crelle's
              Journal]},
    VOLUME = {461},
      YEAR = {1995},
     PAGES = {137--178},
      ISSN = {0075-4102,1435-5345},
   MRCLASS = {46B20 (46B28 47D15)},
   }

@article {C,
    AUTHOR = {Cho, Chong-Man},
     TITLE = {Spaces of compact operators which are {$M$}-ideals in
              {$L(X,Y)$}},
   JOURNAL = {Internat. J. Math. Math. Sci.},
  FJOURNAL = {International Journal of Mathematics and Mathematical
              Sciences},
    VOLUME = {15},
      YEAR = {1992},
    NUMBER = {3},
     PAGES = {617--620},
      ISSN = {0161-1712,1687-0425},
   MRCLASS = {47D15 (46A32 46B28 47B07)},
   }

@book {AK,
    AUTHOR = {Albiac, Fernando and Kalton, Nigel J.},
     TITLE = {Topics in {B}anach space theory},
    SERIES = {Graduate Texts in Mathematics},
    VOLUME = {233},
   EDITION = {Second},
 PUBLISHER = {Springer, [Cham]},
      YEAR = {2016},
     PAGES = {xx+508},
      ISBN = {978-3-319-31555-3; 978-3-319-31557-7},
   MRCLASS = {46-01 (46B04 46B07 46B15 46B20 46B22)},
  MRNUMBER = {3526021},
MRREVIEWER = {Mari\'an\ Fabian},
       DOI = {10.1007/978-3-319-31557-7},
       URL = {https://doi.org/10.1007/978-3-319-31557-7},
}

@article {BS,
    AUTHOR = {Bhatia, Rajendra and \u{S}emrl, Peter},
     TITLE = {Orthogonality of matrices and some distance problems},
   JOURNAL = {Linear Algebra Appl.},
  FJOURNAL = {Linear Algebra and its Applications},
    VOLUME = {287},
      YEAR = {1999},
    NUMBER = {1-3},
     PAGES = {77--85},
      ISSN = {0024-3795,1873-1856},
   MRCLASS = {15A60 (47A30 47B47)}
}

@article {BFS,
    AUTHOR = {Ben\'itez, Carlos and Fern\'andez, Manuel and Soriano, Mar\'ia
              L.},
     TITLE = {Orthogonality of matrices},
   JOURNAL = {Linear Algebra Appl.},
  FJOURNAL = {Linear Algebra and its Applications},
    VOLUME = {422},
      YEAR = {2007},
    NUMBER = {1},
     PAGES = {155--163},
      ISSN = {0024-3795,1873-1856},
   MRCLASS = {46C15 (15A60 47A30)},
}

@article {LS,
    AUTHOR = {Li, Chi-Kwong and Schneider, Hans},
     TITLE = {Orthogonality of matrices},
   JOURNAL = {Linear Algebra Appl.},
  FJOURNAL = {Linear Algebra and its Applications},
    VOLUME = {347},
      YEAR = {2002},
     PAGES = {115--122},
      ISSN = {0024-3795,1873-1856},
   MRCLASS = {15A60 (46B99)},
}

@article {T,
    AUTHOR = {Turn\v{s}ek, Aleksej},
     TITLE = {A remark on orthogonality and symmetry of operators in
              {$\mathcal{B}(\mathcal{H})$}},
   JOURNAL = {Linear Algebra Appl.},
  FJOURNAL = {Linear Algebra and its Applications},
    VOLUME = {535},
      YEAR = {2017},
     PAGES = {141--150},
      ISSN = {0024-3795,1873-1856},
   MRCLASS = {47L05},
}

@article {PSG,
    AUTHOR = {Paul, K. and Sain, D. and Ghosh, P.},
     TITLE = {Birkhoff-{J}ames orthogonality and smoothness of bounded
              linear operators},
   JOURNAL = {Linear Algebra Appl.},
  FJOURNAL = {Linear Algebra and its Applications},
    VOLUME = {506},
      YEAR = {2016},
     PAGES = {551--563},
      ISSN = {0024-3795,1873-1856},
   MRCLASS = {46B20 (46B50)},
}

@article {SPH,
    AUTHOR = {Sain, Debmalya and Paul, Kallol and Hait, Sourav},
     TITLE = {Operator norm attainment and {B}irkhoff-{J}ames orthogonality},
   JOURNAL = {Linear Algebra Appl.},
  FJOURNAL = {Linear Algebra and its Applications},
    VOLUME = {476},
      YEAR = {2015},
     PAGES = {85--97},
      ISSN = {0024-3795,1873-1856},
   MRCLASS = {47A30 (46C15)},
}

@article {MPRS,
    AUTHOR = {Mal, Arpita and Paul, Kallol and Rao, T. S. S. R. K. and Sain,
              Debmalya},
     TITLE = {Approximate {B}irkhoff-{J}ames orthogonality and smoothness in
              the space of bounded linear operators},
   JOURNAL = {Monatsh. Math.},
  FJOURNAL = {Monatshefte f\"ur Mathematik},
    VOLUME = {190},
      YEAR = {2019},
    NUMBER = {3},
     PAGES = {549--558},
      ISSN = {0026-9255,1436-5081},
   MRCLASS = {46B28 (46B20 47L05)},
}

@article {CK,
    AUTHOR = {Choi, Geunsu and Kim, Sun Kwang},
     TITLE = {The {B}irkhoff-{J}ames orthogonality and norm attainment for
              multilinear maps},
   JOURNAL = {J. Math. Anal. Appl.},
  FJOURNAL = {Journal of Mathematical Analysis and Applications},
    VOLUME = {502},
      YEAR = {2021},
    NUMBER = {2},
     PAGES = {Paper No. 125275, 16},
      ISSN = {0022-247X,1096-0813},
   MRCLASS = {46B04 (46B07 46B20 46G25)},
}

@article {FP,
    AUTHOR = {Franchetti, Carlo and Pay\'a, Rafael},
     TITLE = {Banach spaces with strongly subdifferentiable norm},
   JOURNAL = {Boll. Un. Mat. Ital. B (7)},
  FJOURNAL = {Unione Matematica Italiana. Bollettino. B. Serie VII},
    VOLUME = {7},
      YEAR = {1993},
    NUMBER = {1},
     PAGES = {45--70},
MRCLASS = {46B20 (46B22)},
}

@article {B,
    AUTHOR = {Birkhoff, Garrett},
     TITLE = {Orthogonality in linear metric spaces},
   JOURNAL = {Duke Math. J.},
  FJOURNAL = {Duke Mathematical Journal},
    VOLUME = {1},
      YEAR = {1935},
    NUMBER = {2},
     PAGES = {169--172},
      ISSN = {0012-7094,1547-7398},
   MRCLASS = {99-04},
}

@book {HWW,
    AUTHOR = {Harmand, P. and Werner, D. and Werner, W.},
     TITLE = {{$M$}-ideals in {B}anach spaces and {B}anach algebras},
    SERIES = {Lecture Notes in Mathematics},
    VOLUME = {1547},
 PUBLISHER = {Springer-Verlag, Berlin},
      YEAR = {1993},
     PAGES = {viii+387},
      ISBN = {3-540-56814-X},
   MRCLASS = {46Bxx (46-02 46B20 46H99 46J99)},
}

@article{siju, 
title={Differentiability of the operator norm on $\ell _p$ spaces}, 
DOI={}, 
journal={Canadian Journal of Mathematics}, 
author={Siju, Sreejith}, year={2024}, 
pages={1--20}
}

@article {J,
    AUTHOR = {James, Robert C.},
     TITLE = {Orthogonality and linear functionals in normed linear spaces},
   JOURNAL = {Trans. Amer. Math. Soc.},
  FJOURNAL = {Transactions of the American Mathematical Society},
    VOLUME = {61},
      YEAR = {1947},
     PAGES = {265--292},
      ISSN = {0002-9947,1088-6850},
   MRCLASS = {46.0X},
}

@article {KL,
    AUTHOR = {Kim, Sun Kwang and Lee, Han Ju},
     TITLE = {The {B}irkhoff-{J}ames orthogonality of operators on infinite
              dimensional {B}anach spaces},
   JOURNAL = {Linear Algebra Appl.},
  FJOURNAL = {Linear Algebra and its Applications},
    VOLUME = {582},
      YEAR = {2019},
     PAGES = {440--451},
      ISSN = {0024-3795,1873-1856},
   MRCLASS = {46B20 (46B04 46B22 46B28)},
}

@article {K,
    AUTHOR = {Kim, Sun Kwang},
     TITLE = {Quantity of operators with {B}hatia-\v Semrl property},
   JOURNAL = {Linear Algebra Appl.},
  FJOURNAL = {Linear Algebra and its Applications},
    VOLUME = {537},
      YEAR = {2018},
     PAGES = {22--37},
      ISSN = {0024-3795,1873-1856},
   MRCLASS = {46B20 (46B04 46B22)},
}

@article {SPMR,
    AUTHOR = {Sain, Debmalya and Paul, Kallol and Mal, Arpita and Ray,
              Anubhab},
     TITLE = {A complete characterization of smoothness in the space of
              bounded linear operators},
   JOURNAL = {Linear Multilinear Algebra},
  FJOURNAL = {Linear and Multilinear Algebra},
    VOLUME = {68},
      YEAR = {2020},
    NUMBER = {12},
     PAGES = {2484--2494},
      ISSN = {0308-1087,1563-5139},
   MRCLASS = {46B20 (47L05)},
}

@book {BDE,
    AUTHOR = {Behrends, Ehrhard and Danckwerts, Rainer and Evans, Richard
              and G\"obel, Silke and Greim, Peter and Meyfarth, Konrad and
              M\"uller, Winfried},
     TITLE = {{$L\sp{p}$}-structure in real {B}anach spaces},
    SERIES = {Lecture Notes in Mathematics},
    VOLUME = {613},
 PUBLISHER = {Springer-Verlag, Berlin-New York},
      YEAR = {1977},
     PAGES = {x+108},
      ISBN = {3-540-08441-X},
   MRCLASS = {46B99 (46E30)},
}

@article {S,
    AUTHOR = {Singla, Sushil},
     TITLE = {Gateaux derivative of {$C^\ast$} norm},
   JOURNAL = {Linear Algebra Appl.},
  FJOURNAL = {Linear Algebra and its Applications},
    VOLUME = {629},
      YEAR = {2021},
     PAGES = {208--218},
      ISSN = {0024-3795,1873-1856},
   MRCLASS = {46G05 (46B20 46L05 46L30)},
}

@book {I,
    AUTHOR = {Cioranescu, Ioana},
     TITLE = {Geometry of {B}anach spaces, duality mappings and nonlinear
              problems},
    SERIES = {Mathematics and its Applications},
    VOLUME = {62},
 PUBLISHER = {Kluwer Academic Publishers Group, Dordrecht},
      YEAR = {1990},
     PAGES = {xiv+260},
      ISBN = {0-7923-0910-3},
   MRCLASS = {46Bxx (47Hxx 58C30)},
}

@book {M,
    AUTHOR = {Megginson, Robert E.},
     TITLE = {An introduction to {B}anach space theory},
    SERIES = {Graduate Texts in Mathematics},
    VOLUME = {183},
 PUBLISHER = {Springer-Verlag, New York},
      YEAR = {1998},
     PAGES = {xx+596},
      ISBN = {0-387-98431-3},
   MRCLASS = {46-01 (46Bxx 47-01)},
}

@article {MJAP,
    AUTHOR = {Mart\'inez-Moreno, J. and Mena-Jurado, J. F. and
              Pay\'a-Albert, R. and Rodr\'iguez-Palacios, \'A.},
     TITLE = {An approach to numerical ranges without {B}anach algebra
              theory},
   JOURNAL = {Illinois J. Math.},
  FJOURNAL = {Illinois Journal of Mathematics},
    VOLUME = {29},
      YEAR = {1985},
    NUMBER = {4},
     PAGES = {609--626},
      ISSN = {0019-2082,1945-6581},
   MRCLASS = {46H99 (46K99 46L70 47A12)},
}
\bibliographystyle{amsplain}
\end{document}